\documentclass[12pt]{article}

\usepackage{amsfonts,graphicx}

\newcommand{\eh}{\hspace{.06in}}
\newcommand{\cc}{\bold c}

\newcommand{\A}{{\cal{A}}}

\newcommand{\B}{{\cal{B}}}

\newcommand{\F}{{\cal{F}}}

\newcommand{\C}{\mathbb{C}}
\newcommand{\R}{\mathbb{R}}

\newcommand{\hC}{\hat{\C}}
\newcommand{\af}{\alpha}

\newcommand{\ds}{\displaystyle}

\newcommand{\ovl}{\overline}

\newcommand{\BE}{\begin{equation}}
\newcommand{\EE}{\end{equation}}

\newtheorem{thm}{Theorem}[section]
\newtheorem{lem}{Lemma}[section]
\newtheorem{dft}{Definition}[section]
\newtheorem{rmk}{Remark}[section]

\title{Explicit minimal Scherk saddle towers of arbitrary even genera in $\R^3$}
\author{A.J. Yucra Hancco\footnote{Partially supported by CAPES.} $\cdot$ G.A. Lobos $\cdot$ V. Ramos Batista}
\date{}
\begin{document}

\maketitle 
\begin{abstract}
Starting from works by Scherk (1835) and by Enneper-Weierstra\ss \ (1863), new minimal surfaces with Scherk ends were found only in 1988 by Karcher (see \cite{Karcher1,Karcher}). In the singly periodic case, Karcher's examples of positive genera had been unique until Traizet obtained new ones in 1996 (see \cite{Traizet}). However, Traizet's construction is implicit and excludes {\it towers}, namely the desingularisation of more than two concurrent planes. Then, new explicit towers were found only in 2006 by Martin and Ramos Batista (see \cite{Martin}), all of them with genus one. For genus two, the first such towers were constructed in 2010 (see \cite{Valerio2}). Back to 2009, implicit towers of arbitrary genera were found in \cite{HMM}. In our present work we obtain {\it explicit} minimal Scherk saddle towers, for any given genus $2k$, $k\ge3$.
\end{abstract}

\section{Introduction}
Let $S$ be a complete minimal surface embedded in $\R^3$ and of finite total curvature. If $S$ is neither a plane nor a catenoid, the works of Schoen \cite{Schoen} and L\'opez-Ros \cite{Lopez2} show that $S$ must have positive genus and a number of ends $n\le3$. Such an $S$ was unknown until 1984, when Costa obtained his famous example \cite{Costa}. It was later generalised by Hoffman-Meeks and Hoffman-Karcher in \cite{Hoffman1} and \cite{Hoffman}, respectively. Moreover, in \cite{Hoffman1} the authors launched their conjecture that $n\le$ genus+2 for any such $S$, which still remains open after over a quarter of a century.

For $S$ embedded in a flat space, in 1989 Karcher presented several examples that answer many important questions in the Theory of Minimal Surfaces \cite{Karcher1,Karcher}. Among others, he obtained the first $S$ with positive genera and helicoidal ends, proved the existence of Schoen's surfaces \cite{Schoen1} and found singly and doubly periodic $S$ that do not belong to Scherk's minimal surface families. 

By the way, Karcher constructed saddle towers $S$ of genera zero and one (in the quotient by their translation group), and number of ends $n=2k$, $k\ge$ genus+2. Since then, very few new explicit $S$ were found, such as in \cite{Martin,Valerio2}. This can be due to strong restrictions that underlie these surfaces. For instance, Meeks and Wolf proved in \cite{Meeks} that $S$ belongs to Scherk's second family if $n=4$.
\begin{figure}[ht!]
\centering
\includegraphics[scale=0.6]{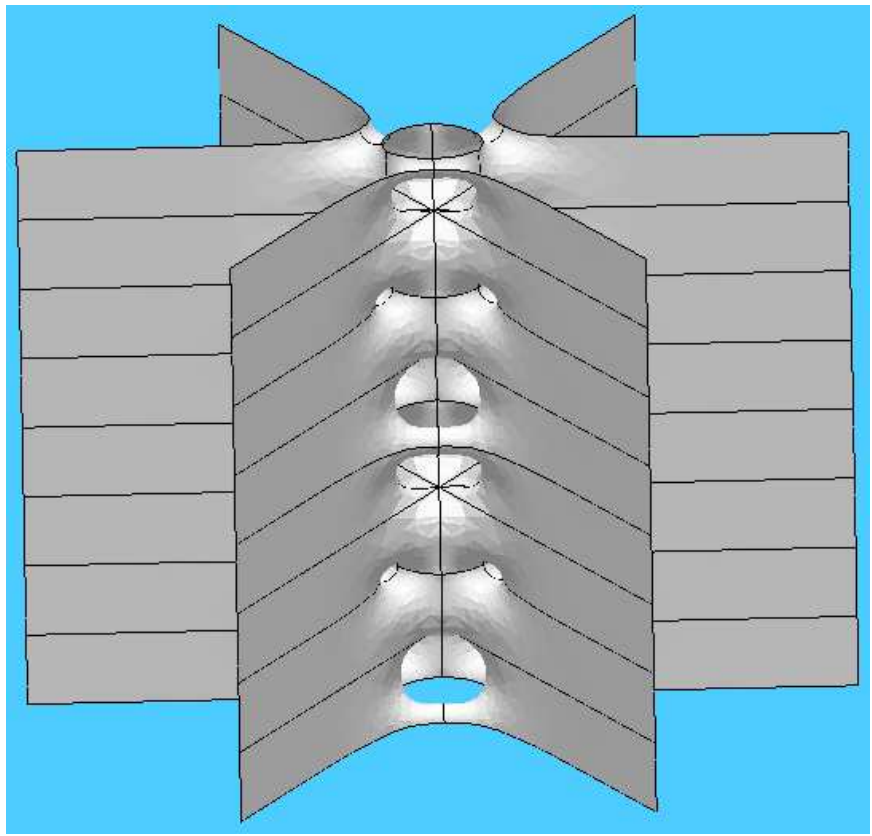}
\caption{A Scherk saddle tower of genus $2k$, $k=3$}.
\end{figure}

In this work we present the first explicit $S$ of arbitrary genera $2k$ and $2k$ Scherk ends, $k\ge3$. More specifically, we prove
\begin{thm}
For each natural $k\ge3$ there exists a continuous one-parameter family of embedded minimal saddle towers in $\R^3$, of which any member $ST_{2k}$ has its symmetry group generated by the following maps:
\begin{enumerate}
\itemsep = 0.0 pc
\parsep  = 0.0 pc
\parskip = 0.0 pc
\item $\pi$-rotation about the line $[(\cot\frac{\pi}{2k},1,0)]\subset\R^3$;
\item Reflection in the vertical plane $Ox_1x_3$;
\item Reflection in the horizontal plane $(0,0,1)+Ox_1x_2$.
\end{enumerate}
Composition of items 1 and 3 make $ST_{2k}$ invariant by the translation group $G=\langle(0,0,4)\rangle$. Moreover, $ST_{2k}/G$ has $2k$ Scherk ends and genus $2k$. The surfaces $ST_{2k}$ are embedded in $\R^3$.
\label{mainthm}
\end{thm}

Notice that items 1 and 2 make $ST_{2k}$ invariant by $\rho$, defined as the composite of $\pi/k$-rotation around $Ox_3$ and reflection in $Ox_1x_2$. The symmetry $\rho$ will be useful in our constructions. 

Regarding explicit saddle towers $S$ with arbitrary odd genus, we are convinced of their existence but we prefer to leave it as an open question, in spite of \cite{HMM}. There the authors constructed implicit examples for any positive genus, for which however the inequality $n\ge$ 2(genus+2) does not hold. Although this could not be another kind of Hoffman-Meeks conjecture, the examples from Theorem 1.1 still verify that inequality like all {\it explicit} examples found to date.

In fact, if one aims at classifying minimal surfaces, then explicit constructions are strictly necessary.

\section{Preliminaries}

This section presents some basic definitions and theorems used throughout this work. We only consider surfaces that are regular and connected. For details see \cite{Conway, Foster, Hoffman, JorgeM, Karcher, Lopez, Nitsche, Osserman}.
 
\begin{thm}Let $X:R\to\mathbb{E}$ be a complete isometric immersion of a Riemann surface $R$ into a three-dimensional flat space $\mathbb{E}$. If $X$ is minimal and the total Gaussian curvature $\int_RKdA$ is finite, then $R$ is conformal to $\ovl{R}\setminus\{p_1,p_2,\cdots,p_r\}$, where $\ovl{R}$ is a compact Riemann surface and $r$ is a certain number of points $\{p_1,p_2,\cdots,p_r\}\subset\ovl{R}$.\label{hueber}
\end{thm}

\begin{thm}{\rm (Weierstra\ss \ Representation)}
Let $R$ be a Riemann surface, $g$ and $dh$ meromorphic function and 1-differential form on $R$, such that the zeros of $dh$ coincide with the poles and zeros of $g$. Suppose that $X:R\to\mathbb{E}$ given by
\BE
   X(p):={\rm Re}\int^p\Phi,\eh\eh\Phi=\frac{1}{2}(1/g-g,i/g+ig,2)dh,
   \label{eqw}
\EE
is well-defined. The $X$ is a conformal minimal immersion. Conversely, every conformal minimal immersion $X:R\to\mathbb{E}$ can be expressed as $(\ref{eqw})$ for some meromorphic function $g$ and 1-form $dh$.\label{WR}
\end{thm}

\begin{dft}\rm The pair $(g,dh)$ is the {\it Weierstra\ss \ data} and the components of $\Phi=(\phi_1,\phi_2,\phi_3)$ are the {\it Weierstra\ss \ forms} on $R$ of the minimal immersion $X:R\to X(R)\subset\mathbb{E}$.\label{names2WR} 
\end{dft}

\begin{dft}\rm Let $R$ and $\ovl{R}$ be as in Theorems \ref{hueber} and \ref{WR}. An {\it end} of $R$ is the image by $X$ of a punched neighbourhood $V_p$, $p\in\{p_1,p_2,\cdots,p_r\}$, such that $(\{p_1,p_2,\cdots,p_r\}\setminus\{p\})\cap\ovl{V}_p=\emptyset$. The end is embedded if $X:V_p\to\mathbb{E}$ is an embedding for a sufficiently small $V_p$.\label{endofms} 
\end{dft}

\begin{thm}Under the hypotheses of Theorems \ref{hueber} and \ref{WR}, the Weierstra\ss \ data $(g,dh)$ extend meromorphically on $\ovl{R}$.\label{isalg}
\end{thm}

\begin{thm}Let $\mathbb{X}$ and $\mathbb{Y}$ be Riemann surfaces and $f:\mathbb{X}\to\mathbb{Y}$ a non-constant proper holomorphic map. In this case, there is a natural number $n$ such that $f$ attains each point $q\in\mathbb{Y}$ exactly $n$ times, including multiplicity.
\end{thm}

\begin{dft}\rm Let $\mathbb{X}$ and $\mathbb{Y}$ be Riemann surfaces and $f:\mathbb{X}\to\mathbb{Y}$ a non-constant meromorphic function. The {\it degree} of $f$ is the cardinality of $f^{-1}(q)$, $\forall\,q\in\mathbb{Y}$, denoted by {\it deg}$(f)$.
\end{dft}

\begin{thm}{\rm (Jorge-Meeks Formula)} Let $X:R\to\mathbb{E}$ be a complete minimal surface with finite total curvature $\int_RKdA$. If $R$ has ends that are all embedded, then $deg(g)=k+r-1$, where $k$ is the genus of $\ovl{R}=R\cup\{p_1,p_2,\cdots,p_r\}$ and $r$ is the number of ends.\label{JM}
\end{thm}

\begin{rmk}\rm In the proof of Theorem \ref{JM}, for the case of Scherk-ends the variable $r$ counts them in pairs. The function $g$ is the stereographic projection of the Gau\ss \ map $N:R\to S^2$ of the minimal immersion $X$. It is a branched covering map of $\hC$ and $\int_RKdA=-4\pi deg(g)$.\label{remark}
\end{rmk}

\begin{thm}If $\sigma$ is a curve on $X(R)$, then
\begin{enumerate}
\itemsep = 0.0 pc
\parsep  = 0.0 pc
\parskip = 0.0 pc
\item [i)] $\sigma$ is asymptotic if and only if $(dh\cdot dg/g)|_{\sigma'}\in i\R$;
\item[ii)] $\sigma$ is a principal curvature line if and only if $(dh\cdot dg/g)|_{\sigma'}\in\R$.
\end{enumerate}
\label{dhdgg}
\end{thm}

\begin{thm}{\rm (Schwarz Reflection Principle)} Let $S$ be a complete minimal surface. If $\ell$ is a straight line on $S$, then $S$ is invariant by 180$^\circ$-rotation about $\ell$. If $\af$ is a planar geodesic on $S$, then $S$ is invariant by reflection in the plane of $\af$.\label{schwarz}
\end{thm} 

\begin{thm}If in some holomorphic coordinates of a minimal immersion $F:\Omega\to\R^3$ there is a curve $\sigma$ such that $g(\sigma)$ is either in a meridian or in the equator of $\hC=S^2$, and $dh(\sigma')\subset\R\cup i\R$, then $F\circ\sigma$ is either in a plane or in a straight line. In both cases, $\sigma$ is a geodesic. It is planar exactly when $(dh\cdot dg/g)|_{\sigma'}\in\R$ and it is straight exactly when $(dh\cdot dg/g)|_{\sigma'}\in i\R$.\label{meridequator}
\end{thm}

\begin{thm}
Let $G$ be a open connected subset of $\C$. Then $G$ is simply connected if and only if $\hC\setminus G$ is connected.\label{simcon}
\end{thm}

\section{Construction of the surfaces $ST_{2k}$}

In Theorem \ref{mainthm} we denoted our surfaces by $ST_{2k}$. This theorem is proved by Karcher's reverse construction method. Namely, we derive a list of necessary conditions that must hold in case the surfaces exist. They will end up in algebraic equations for $R$, $g$ and $dh$. At this point, Theorems \ref{hueber} and \ref{WR} apply. Afterwards, we must prove that $X:R\to\mathbb{E}$ really corresponds to each surface $ST_{2k}$ from Theorem \ref{mainthm}. 

Suppose we had a minimal surface like in Figure 1. Take the quotient by its translation group, followed by a compactification of the ends. We get a fundamental piece $\ovl{S}$. 

Now $\ovl{S}$ has genus $2k$, and we assume that $\ovl{S}$ is invariant by $\pi/k$-rotation around $Ox_3$ {\it followed} by a reflection in $Ox_1x_2$. Let us  denote this symmetry by $\rho$. Hence, the Euler characteristic of $\rho(\ovl{S})$ is
\[
   \chi(\rho(\ovl{S}))=\frac{\chi(\ovl{S})}{2k}+\frac{2k-1}{2k}\cdot2=\frac{1}{k}-2+2-\frac{1}{k}=0.
\]
Therefore, $\rho(\ovl{S})$ is a torus $T$. Due to the horizontal reflectional symmetries of $\ovl{S}$, $T$ is a rectangular torus. 

\subsection{The function $z$ on $\ovl{S}$ and the Gau\ss \ map $g(z)$}

We have just obtained the rectangular torus $T=\rho(\ovl{S})$. Let us now obtain two meromorphic functions $g$, $z$ on $\ovl{S}$ through the pullback by $\rho$ of functions on $T$. Because of Remark \ref{remark}, $g$ will be constructed by looking at the stereographic projection of the Gau\ss \ map $N:\ovl{S}\to S^2=\hC$. Regarding $z$, we choose it to make the relation $g=g(z)$ as simple as possible. See Figure 2 for an illustration.

\begin{figure}[h!]
\centering
\includegraphics[scale=0.80]{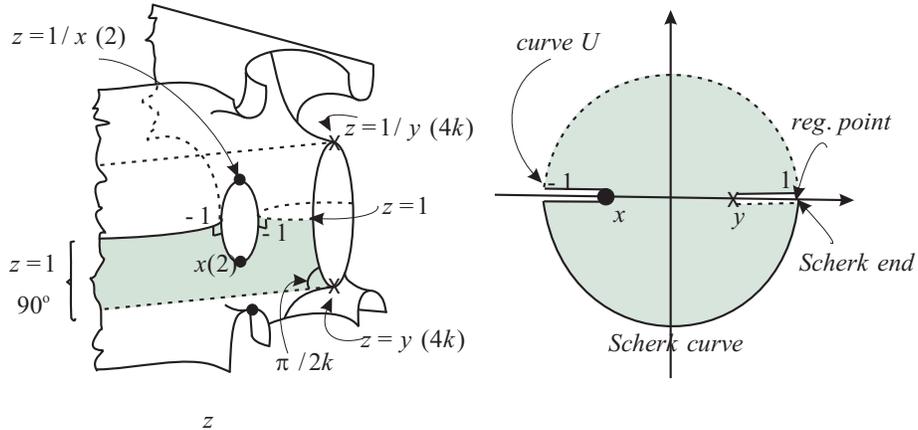}
\caption{Values of $z$ at special points of $\ovl{S}$.}
\end{figure}
 
In \cite{Valerio1} the author considers an elliptic function ${Z'}:T\to\C$ as schematised in Figure 3. We define $z:\ovl{S}\to\hC$ as $z:=Z'\circ\rho$.

\begin{figure}[h!]
\hspace{1cm}
\centering
\includegraphics[scale=0.80]{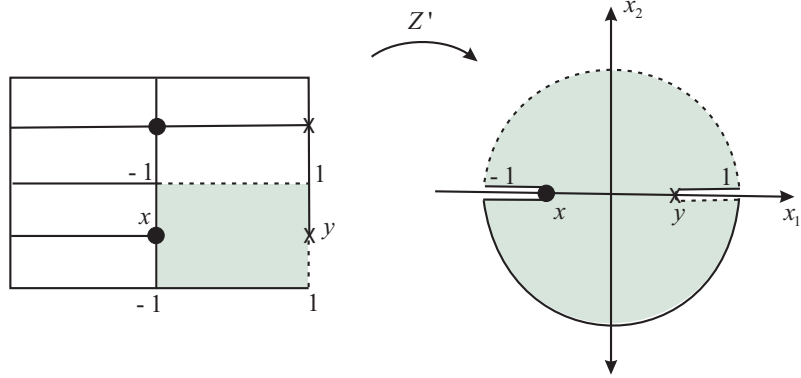}
\caption{The torus $T$ with some special values of $Z'$.}
\end{figure}

By looking at the normal vector on $\ovl{S}$, we know that $g$ has poles and zeros as shown in Figure 4.

\begin{figure}[h!]
\centering
\includegraphics[scale=0.80]{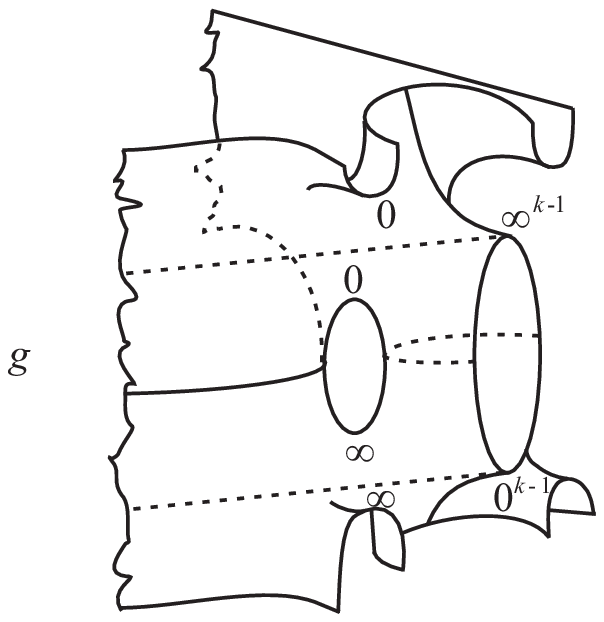}
\caption{Poles and zeros of $g$ on $\ovl{S}$.}
\end{figure}
 
The corresponding values of $z$ and $g$ are represented in Figures 2 and 4, including multiplicities. Therefore, we obtain the following algebraic relation between $g$ and $z$:
\BE
     g^{4k}=c\left(\frac{z-y}{1-yz}\right)^{k-1}\left(\frac{1-xz}{x-z}\right)^{2k}\label{fg},
\EE
where $c$ is a real constant. From Figure 2, we see that along $z(t)=e^{i\pi t}$, $0<t<1$, the following holds: $|g|=1\Longleftrightarrow|z|=1$. Therefore,
\[
   |g^{4k}|=\biggl|c\left(\frac{y-z}{1-yz}\right)^{k-1}\left(\frac{1-xz}{x-z}\right)^{2k}\biggl|,
\]
which implies $|c|=1$. For $z=1$ we have $g=e^{i\pi(1+k)/k}$, hence $c=1$. 

\begin{figure}[h!]
\centering
\includegraphics[scale=0.80]{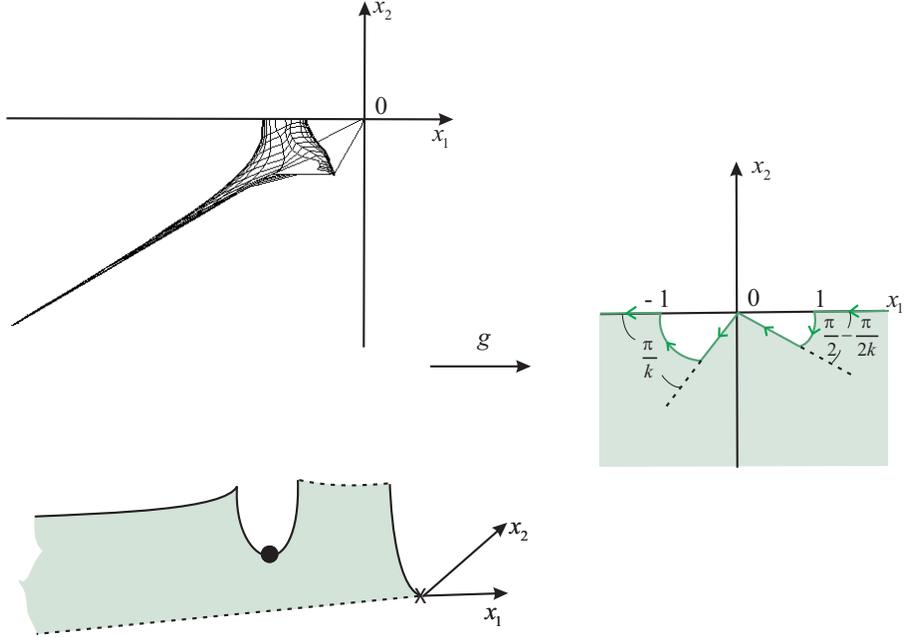}\\
\caption{The fundamental domain viewed from above and the corresponding values of $g$.}
\end{figure}
 
From (\ref{fg}), and based on Figures 2 and 5, we summarise what we know about $g$ and $dh$ in Table \ref{gsymm}. Of course, we have not defined $dh$ yet but it is determined by $g$, according to the theorems listed in Section 2.

\begin{table}[ht]
\centering
\begin{tabular}{|c||c|c|c|}
\hline
{\rm stretch} & {\rm $z$-values} &{\rm $g$-values} & {\rm $dh(z')$-values} \\\hline\hline
1             & $z(t)=e^{i\pi t},0<t<1$  & $|g|=1$    & $\in i\R$             \\\hline 
2             & $-1<z<x$         & $g\leq -1$      & $\notin\R\cup i\R$    \\\hline
3             & $x<z<-1$         & $g\geq 1$       & $\notin\R\cup i\R$    \\\hline
4             & $z(t)=e^{i\pi t},-1<t<0$ & $|g|=1$    & $\in i\R$             \\\hline
5             & $1<z<y $         & $\eh\eh \in -ie^{.5\pi i/k}\R$ & $\in i\R$  \\\hline
6             & $y<z<1 $         & $\eh \in-e^{\pi i/k}\R$ & $\in\R$          \\\hline
\end{tabular}
\caption{Values of $g$ and $dh(z')$ along the symmetry curves.}
\label{gsymm}
\end{table}

\begin{rmk}\rm As we shall see in Section \ref{secdh}, along the stretches 2 and 3 from Table \ref{gsymm}, $dh(z')$ takes complex values not in $\R\cup i\R$. In fact, we want to apply Theorems \ref{dhdgg} and \ref{meridequator}, which give sufficient conditions to prove the existence of symmetry curves. However, $ST_{2k}$ will not have more symmetries than the ones listed in Theorem \ref{mainthm}.\label{obsU}
\end{rmk}

Table \ref{gsymm} shows that $g$ is consistent with the normal vector on $\ovl{S}$ along special curves on the surface. Now we list important involutions of $\ovl{S}$ related with the symmetries of $ST_{2k}$ in $\R^3$. Based on (\ref{fg}) and Figures 2, 5 we can summarise these involutions in Table \ref{inv}:
 
\begin{table}[ht]
\centering
\begin{tabular}{|c||c|c|}
\hline
{\rm symmetry}        & {\rm involution}                     & $g\in$         \\\hline\hline
$z(t)=e^{i\pi t},0<t<1$  & $(z,g)\to(1/\bar{z},1/\bar{g})$      & $S^1$          \\\hline
$-1<z<x$              & $(z,g)\to(\bar{z},\bar{g})$          & $(-\infty,-1)$ \\\hline
$x<z<-1$              & $(z,g)\to(\bar{z},\bar{g})$          & $(1,\infty) $  \\\hline
$z(t)=e^{i\pi t},-1<t<0$ & $(z,g)\to(1/\bar{z},1/\bar{g})$      & $S^1$          \\\hline
$1<z<y$               &  $(z,g)\to(\bar{z},-e^{i\pi/k}\bar{g})$ & $-ie^{.5i\pi/k}\R$\\\hline
$y<z<1$               & $(z,g)\to(\bar{z},e^{2i\pi/k}\bar{g})$  & $-e^{i\pi/k}\R$   \\\hline
\end{tabular}
\caption{Involutions on $\ovl{S}$.}
\label{inv}
\end{table}

\section{The differential $dh$ in terms of $z$}
\label{secdh}

Since $ST_{2k}$ has Scherk-ends, their corresponding points of $\ovl{S}$ are exactly the poles of $dh$. Regarding the zeros of $dh$, they coincide with the points of $\ovl{S}$ at which $g=0$ or $g=\infty$, including multiplicity. We shall have to read off information about $dz$ in order to write down an equation for $dh$.

\begin{figure}[ht]
\centering
\includegraphics[scale=0.80]{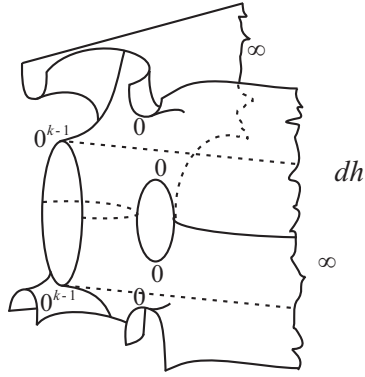}
\caption{Poles and zeros of $dh$ on $\ovl{S}$.}
\end{figure}

Figure 2 shows two points marked with $\times$ at which $z=y$ and $z=1/y$. There we have $dz=0$ of order $4k-1$. Moreover, at $z=x$ we have $dz=0$ of order 1. Now Figure 7 illustrates the divisor of $dh$.

\begin{figure}[ht]
\centering
\includegraphics[scale=0.80]{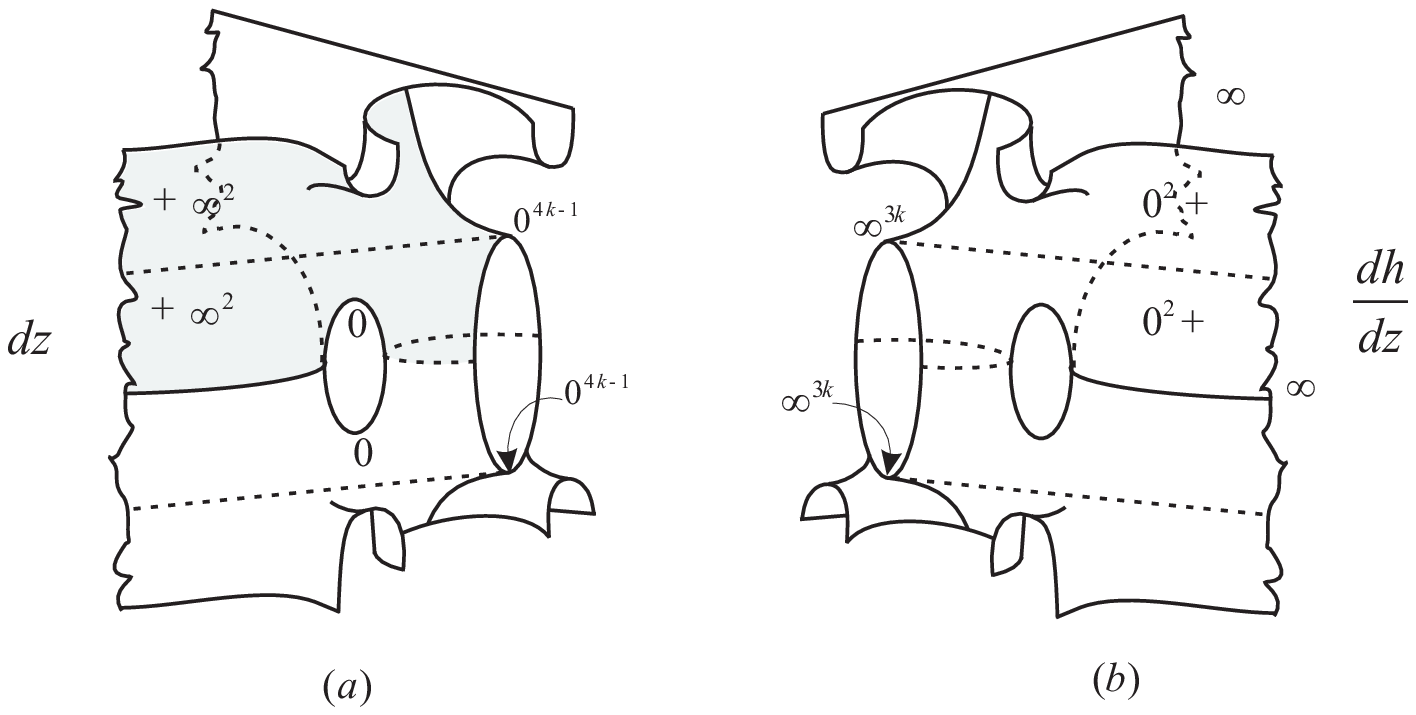}
\caption{(a) Divisor of $dz$ on $\ovl{S}$; (b) divisor of $dh/dz$ on $\ovl{S}$.} 
\end{figure}

In order to obtain $dh$ by means of $dz$, we must analyse the divisor of $f:=(1-yz)(y-z)$. According to Figure 7(a), it is sufficient to construct $F:\ovl{S}\to\hC$ such that $F=f\cdot dh/dz$. Once we have $F$, then
\[
   dh=Fdz/f.
\]

Notice that there are distinct points of $\ovl{S}$ at which $z$ takes the same value 1, namely at the ends and at certain regular points. Since $dh$ has no poles except for the ends, then we shall have to introduce new functions besides $z$ and $g$. They are depicted in Figure 9. In this figure
\[
   v_2=1-\frac{i}{v_1}\biggl(\frac{1-y}{1+y}\biggl),
\]
where
\[
   v_1:=\sqrt{\biggl(\frac{1-x}{1+x}\biggl)^{2}-\biggl(\frac{1-y}{1+y}\biggl)^2}.
\]

\eject
\begin{figure}[h!]
\centering
\includegraphics[scale=0.75]{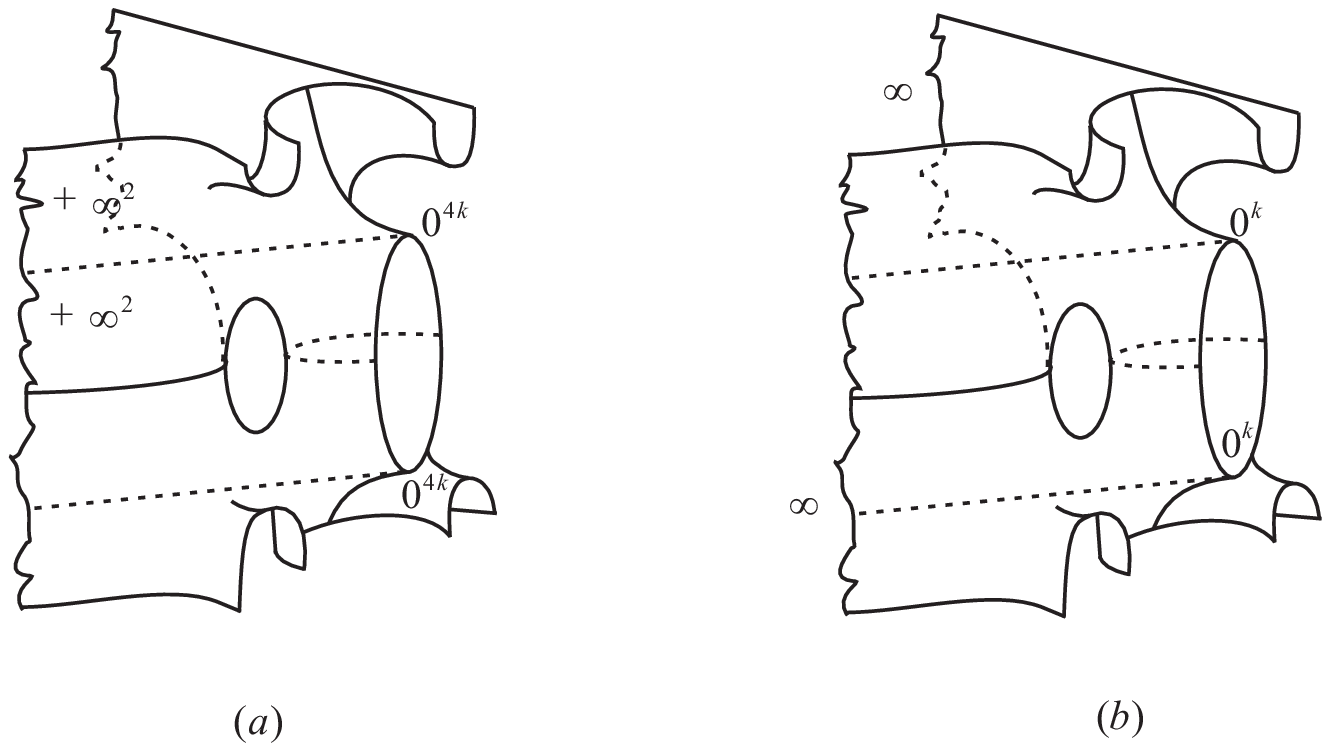}
\caption{(a) Values of $f$; \eh (b) values of $f\cdot dh/dz$ on $\ovl{S}$.}
\end{figure}

\begin{figure}[h!]
\centering
\includegraphics[scale=0.74]{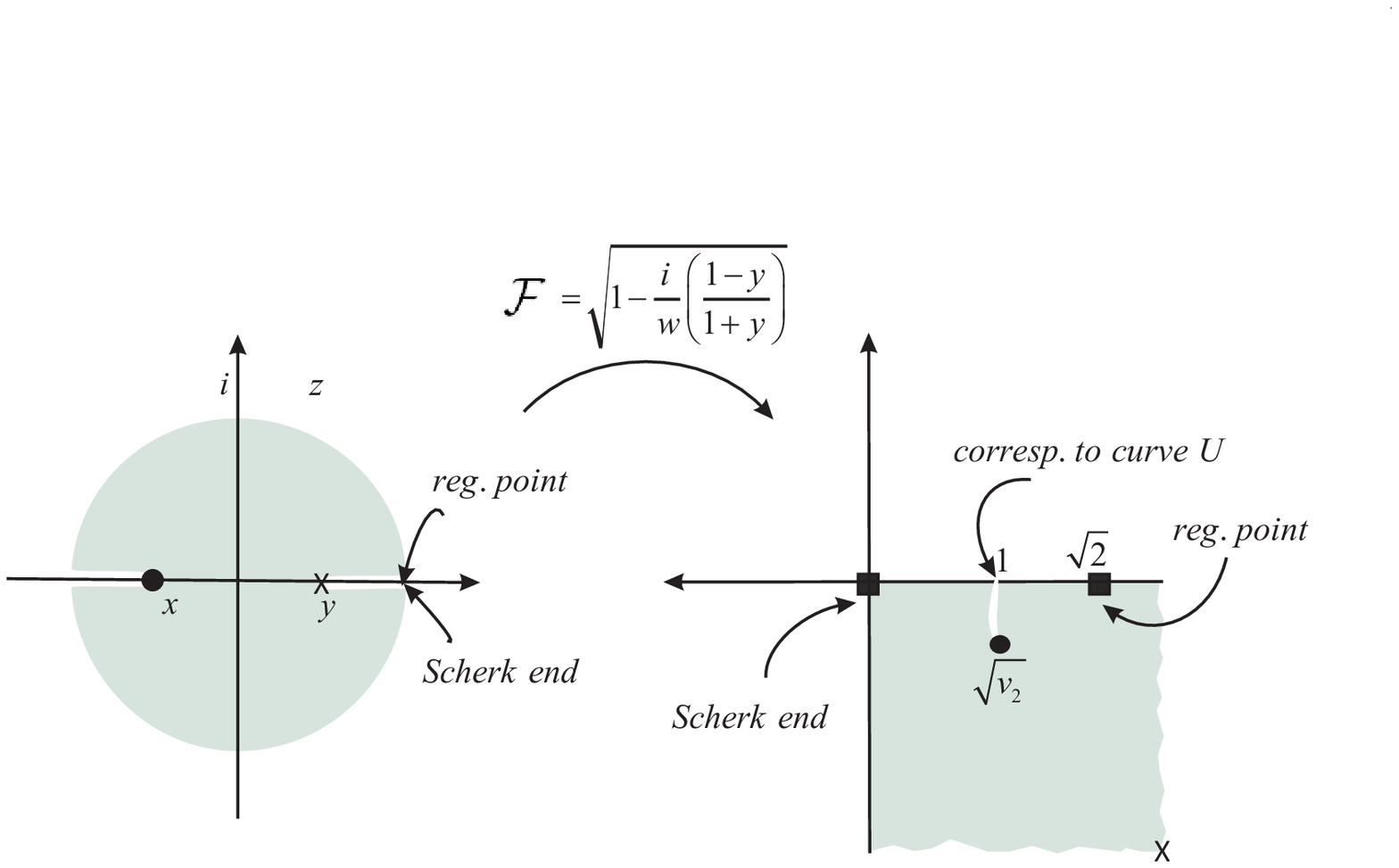}
\caption{The function $\F=\sqrt{1-\frac{i}{w}\left(\frac{1-y}{1+y}\right)}$.}
\end{figure}
\eject

Now
\[
   w:=\sqrt{\biggl(\frac{1-z}{1+z}\biggl)^{2}-\biggl(\frac{1-y}{1+y}\biggl)^2},
\]
is well-defined in a branched covering of $\ovl{S}$ that we call $\ovl{R}$. The algebraic equation of $\ovl{R}$ will be discussed later on. This way we get
\[
   \F=\sqrt{1-\frac{i}{w}\biggl(\frac{1-y}{1+y}\biggl)}.
\]

Therefore, $\F\cdot f\cdot dh/dz$ has neither poles nor zeros, whence must be a non-zero complex constant $\cc$. Namely,
\BE
   dh=\cc\cdot Fdz/f,\label{eqdh}
\EE
where $F:=1/\F$. Now we show that $\cc=1$. Indeed, since stretch 5 in Table 1 is represented by a straight line, the 3rd coordinate of (\ref{eqw}) must be zero. But along this stretch we have $Fdz/f\in i\R$, and in (\ref{eqw}) we compute the real part of a complex integral. Hence $\cc$ must be real. The property of a surface being minimal in $\R^3$ is preserved by the antipodal map and by homotheties. Therefore we can take $\cc=1$.
  
\begin{figure}[ht!]
\centering
\includegraphics[scale=0.90]{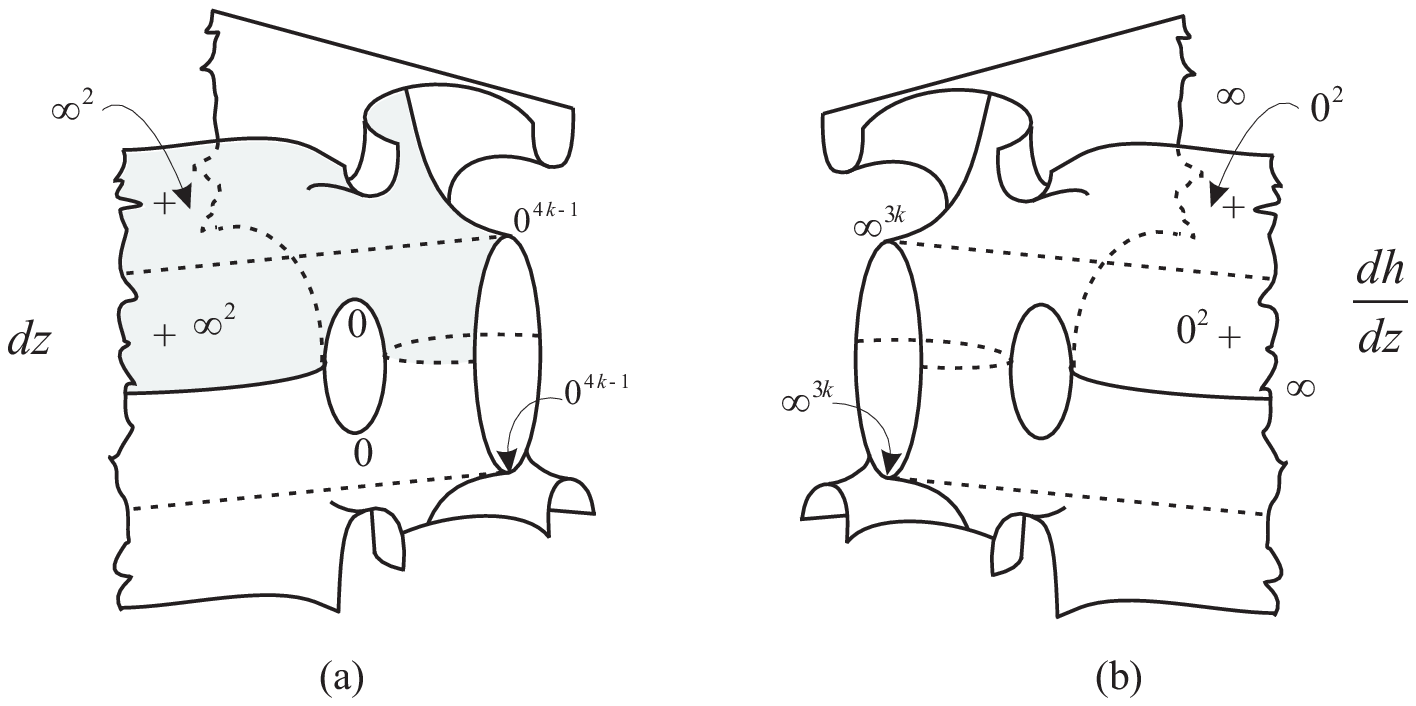}
\caption{Values of $\F$ on $\ovl{S}$.}
\end{figure}

According to Theorems \ref{hueber} and \ref{WR}, $dh$ can be defined on $\ovl{S}$ by a rational function involving only $g,dg,z$ and $dz$. However, its formula is probably far too extensive. Then we use square roots as explained above. Of course, they are not well-defined on $\ovl{S}$, but on a branched covering that we call $\ovl{R}$. In order to describe $\ovl{R}$ by an algebraic equation, we consider
\[
   \F=\sqrt{1-\frac{i}{w}\biggl(\frac{1-y}{1+y}\biggl)}
\]
and
\[
   w^2=\frac{4f}{(1+z)^2(1+y)^2},
\]
whence
\[
   (\F^2-1)^2=-\frac{1}{4}\frac{(1+z)^2(1-y)^2}{(y-z)(1-yz)}.
\]

This way we get the polynomial
\[
   az^2+bz+c=0,
\]
where $a$, $b$ and $c$ depend on $\cal{F}$ and on some complex constants. This results in
\BE
   z=\frac{-b\pm\sqrt{\Delta}}{2a}.\label{bask}
\EE

Now (\ref{fg}) can be rewritten as
\BE
   A_{3k-1}(g)\cdot z^{3k-1}+\ldots+A_0(g)=0,\label{almostR}
\EE
where $A_j$ is a polynomial in $g$, $\forall\,j$. By applying (\ref{bask}) to (\ref{almostR}) we get
\BE
   \pm\sqrt{\Delta}\cdot E_1=E_2,\label{takesq}
\EE
where $E_1$ and $E_2$ are polynomials in $g$ and $\cal{F}$. We square both sides of (\ref{takesq}) and finally get a polynomial $P(g,\F)=0$, which gives an algebraic equation for $\ovl{R}$. The functions $g$ and $\cal{F}$ are then well-defined on $\ovl{R}$. Of course, there is a projection $B:\ovl{R}\to\ovl{S}$ given by $B(g,\F)=g$. However, there is no projection that makes $\cal{F}$ well-defined on $\ovl{S}$, since we must use square roots to equate $\cal{F}$ on $\ovl{S}$.

We recall (\ref{eqdh}) and see that $F/f$ is a rational function on $\ovl{R}$. However, all computations there would have to match the computations on $\ovl{S}$ that use square roots, because $\ovl{R}$ {\it was} obtained from them. Therefore, we shall keep on working with the square roots.

Now we can analyse $dh$ along the symmetry curves in Table 1. Observe that
\BE
   dh=\frac{z}{f}\cdot\frac{1}\F\cdot\frac{dz}{z}.\label{dhexp}
\EE

On the stretch $y<z<1$ the function $\cal{F}$ is real and positive, whereas $f$ is real and negative. Since the curve is $z(t)=t$, we have $dh(z')\in\R$. For $1<z<y$, $-i\cal{F}$ is real and negative, thus $dh(z')\in i\R$. 

Regarding stretches 1 and 4, there we have $z(t)=e^{i\pi t}$, $0<t<1$ and $-1<t<0$, respectively. Notice that $z/f=y^{-1}/(1/z+z-(1/y+y))$ and on these stretches $z/f$ and $\cal{F}$ are both real. Therefore, $dh(z')\in i\R$ because $dz/z\in i\R$.

From Table 1 we see that the $z$-curves 2 and 3 are geodesics according to Theorem \ref{meridequator}. Moreover, the geodesics are planar curves in cases 1, 4 and 6 because $(dh\cdot dg/g)|_{z'}\in\R$, and a straight line in case 5 because $(dh\cdot dg/g)|_{z'}\in i\R$. Therefore, our minimal surfaces $ST_{2k}$ are symmetric with respect to 1, 4, 5 and 6.

But we recall Remark \ref{obsU} regarding curves 2 and 3. From Table 1, $dh(z')\notin\R\cup i\R$. Indeed, for $z(t)=t$, $-1<t<x$, we have $z/f\in\R$ and $\F\notin\R\cup i\R$. Hence $dh(z')\notin\R\cup i\R$. The surfaces $ST_{2k}$ are {\it not} symmetric with respect to the curves 2 and 3, as we shall prove later.

\section{The period problem}

Figure 5 shows the fundamental domain of $ST_{2k}$. Some important details are reproduced again in Figure 11(a), but there we indicate a path that begins at a point marked with $\times$. The path goes upwards and then from the right to the left-hand side, where we find its end. If the $\times$-point is the origin of $\R^3$, then the end of the path ought to be in the plane $x_2=0$. This path is what we call a {\it period curve}.

\begin{figure}[ht!]
\centering
\includegraphics[scale=0.90]{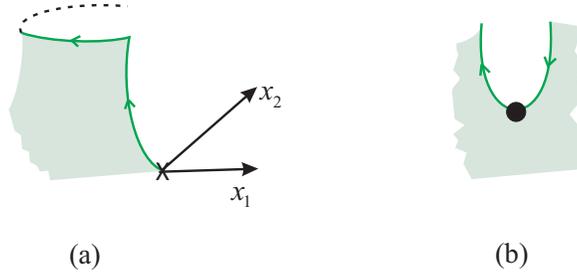}
\caption{Period curves.}
\end{figure}
 
Another one is the U-curve in Figure 11(b). Its period is zero if it has both extremes at the same height. This is determined by the integral of $dh$ alone, and $dh$ is given by square roots. What happens is that their signs change at the vertex of the U-curve, marked with a bullet in Figure 11. The change of sign automatically implies that both extremes of the U-curve {\it do} attain the same height. 

We recall that $dh$ has an {\it algebraic} expression in $\ovl{S}$ involving $(z,g,dz,dg)$. This is ensured by Theorems \ref{hueber} and \ref{WR}. Since $dh=dh(z,g,dz,dg)$ on $\ovl{S}$ is the general expression, it involves extra complex parameters. For the U-curve, its extremes will have heights that depend on these parameters, and the heights will not coincide in general. However, by using square roots in the local expression of $dh$ on $\ovl{S}$, the extra parameters are forced to assume constant values. In fact, we do not even discuss them because of the straight choice of (\ref{eqdh}). But there still remain the two free parameters $x$ and $y$, which also take part in the algebraic equation of $\ovl{R}$.

Hence, we only have {\it one} period problem. For convenience of the reader, herewith we reproduce the Weierstra\ss~data:
\BE
   g^{4k}=(-1)^{k-1}\biggl(\frac{y-z}{1-yz}\biggl)^{k-1}\biggl(\frac{1-xz}{x-z}\biggl)^{2k},
   \label{forgg}
\EE
\BE
   dh=\frac{1}{\sqrt{1-\frac{i(1-y)}{w(1+y)}}}\frac{dz}{(1-yz)(y-z)}, 
   \label{fordh}
\EE
where
\BE
   w=\sqrt{\biggl(\frac{1-z}{1+z}\biggl)^2-\biggl(\frac{1-y}{1+y}\biggl)^2}.\label{forw}
\EE

Let us now analyse some special stretches depicted in Figure 11(a). According to the branches of square root that we have chosen, $1=e^{4(k+1)\pi i}$ for $z(t)=t,\eh y\leq t\leq 1$, which is the upward stretch from $\times$ to a point in $z^{-1}(1)$. Hence
\BE
   g=-e^{i\pi/k}\biggl(\frac{t-y}{1-yt}\biggl)^{\frac{k-1}{4k}}\biggl(\frac{1-xt}{t-x}\biggl)^{\frac{1}{2}}.
   \label{forg1}
\EE

Along the stretch $z(t)=e^{it}$, $0\leq t\leq\pi$, we have 
\BE
   g=-e^{\frac{i\pi}{k}}\cdot e^{\frac{i(k+1)t}{4k}}\biggl(\frac{y-e^{it}}{y-e^{-it}}\biggl)^{\frac{k-1}{4k}}
   \biggl(\frac{x-e^{-it}}{x-e^{it}}\biggl)^{\frac{1}{2}}.\label{forg3}
\EE

Now we analyse $w$ more carefully along $z(t)=t$, $y\leq t\leq1$. From (\ref{forw}) we have
\[
   w=-\frac{2i\sqrt{(1-yt)(t-y)}}{(1+t)(1+y)}.
\]
Let us define $Y:=i(1-y)/(1+y)$. Hence
\[
   \frac{Y}{w}=-\frac{1}{2}\frac{(1+t)(1-y)}{\sqrt{(1-yt)(t-y)}},
\]
and therefore
\[
   dh=\frac{1}{\sqrt{1+\frac{(1+t)(1-y)/2}{\sqrt{(1-yt)(t-y)}}}}\frac{dt}{(1-yt)(y-t)}.
\]
Finally,
\BE
   dh=\frac{1}{\{[(1-yt)(t-y)]^{1/2}+(1+t)(1-y)/2\}^{1/2}}\cdot\frac{dt}{[(1-yt)(t-y)]^{3/4}}.
   \label{dh1}
\EE

In general, (\ref{forw}) rewrites as 
\[
   \frac{w}{Y}=-\frac{2i\sqrt{(1-yz)(z-y)}}{(1+z)(1+y)}\frac{(1+y)}{i(1-y)}=
  -\frac{2\sqrt{y}\sqrt{(1/y+y)-(1/z+z)}}{\sqrt{z}(1/z+1)(1-y)},
\]
and
\[
\left.
\begin{array}{ccl}
dh &=& 
\frac{\ds 1/y}{\sqrt{\ds 1-\frac{\ds i(1-y)}{\ds w(1+y)}}}\frac{\ds dz/z}{\ds (1/y-z)(y/z-1)}\\
   &=&
\frac{\ds 1/y}{\sqrt{\ds 1-\frac{\ds i(1-y)}{\ds w(1+y)}}}\frac{\ds dz/z}{\ds (1/z+z)-(1/y+y)}.\\
\end{array}
\right.
\]

For $z(t)=e^{it}$, $0\leq t\leq \pi$, this results in 
\[
   \frac{w}{Y}=-\frac{\sqrt{y}\sqrt{-2\cos t+(1/y+y)}}{(1-y)\cos (t/2)}\Rightarrow
\]
\[
   \sqrt{1-\frac{Y}{w}}=\sqrt{1+\frac{(1-y)(\cos (t/2))}{\sqrt{y}(1/y+y-2\cos t)^{1/2}}}.
\]

Hence
\BE 
   dh=\frac{i/y}{\sqrt{1+\frac{(1-y)\cos(t/2)}{\sqrt{y}(1/y+y-2\cos t)^{1/2}}}}
      \frac{dt}{2\cos t-(1/y+y)}.\label{fordh3}
\EE

Associated with $z(t)=t,\eh y\leq t\leq 1$, we have
\BE 
   I_1^k=Re\int_y^1\phi_2dh=Re\int_y^1\biggl(\frac{i}{g}+ig\biggl)dh.\label{fori1}
\EE

From (\ref{forg1}), (\ref{dh1}) and the change of variables $t=y+s^{4k}$ we have
\[
\left.I_1^k\right|_{(y,x)}=Re\eh i\int^{(1-y)^{1/4k}}_{0}\left\{\left[-e^{-i\pi/k}\left(\frac{1}{1-y^2-s^{4k}y}\right)^{-\frac{k-1}{4k}}\left(\frac{1-x(y+s^{4k})}{(y+s^{4k})-x}\right)^{-1/2}\right.\right.
\]
\[
\left.\left.-e^{i\pi/k}s^{2(k-1)} \left(\frac{1}{1-y^2-s^{4k}y}\right)^{\frac{k-1}{4k}}\left(\frac{1-x(y+s^{4k})}{(y+s^{4k})-x}\right)^{1/2}\right]\right.\]
\[
\left.\frac{1}{\sqrt{\sqrt{(1-y^2-s^{4k}y)(s^{4k})}+\frac{(1+y+s^{4k})(1-y)}{2}}}\frac{4kds}{\sqrt[4]{[(1-y^2-s^{4k}y)]^{3}}}\right\}.
\]

We recall that $y\in(0,1)$ and $x\in(-1,0)$. For $(y,x)\to(0,-1)$ we have
\BE
\left.
I_1^k\right|_{(0,-1)}=
Re\left\{\eh 4k\sqrt{2}i\int^1_0\left[-e^{-i\pi/k}-e^{i\pi/k} s^{2(k-1)}\right]\frac{ds}{1+s^{2k}}
\right\}.\label{genIk1}
\EE

Let us observe what happens to (\ref{genIk1}) when $k=3$:
\[\left.
I_1^k\right|_{(0,-1)}=Re\left\{\eh 12\sqrt{2}i\left[-e^{-i\pi/3}
\int^1_0\frac{ds}{1+s^6}-e^{i\pi/3}\int^1_0\frac{s^4ds}{1+s^6}\right]\right\}.
\]

Since
\[
\int\frac{ds}{1+s^6}=
\frac{1}{12}\left(-\sqrt{3}\ln{(s^2-\sqrt{3}s+1)}+\sqrt{3}\ln{(s^2+\sqrt{3}s+1)}\right.
\] 
\[
\left.
-2\arctan(\sqrt{3}-2s)+4\arctan(s)+2\arctan(2s+\sqrt{3})\right),
\]
\[
\int\frac{s^4ds}{1+s^6}=\frac{1}{12}\left(\sqrt{3}\ln{(s^2-\sqrt{3}s+1)}-\sqrt{3}\ln{(s^2+\sqrt{3}s+1)}\right.
\]
\[
\left.
-2\arctan(\sqrt{3}-2s)+4\arctan(s)+2\arctan(2s+\sqrt{3})\right),
\] 
then 
\[
\left.
I_1^k\right|_{(0,-1)}\approx 
Re\left\{\eh i12\sqrt{2}\left[e^{-i\pi/3}(0.90377)+e^{i\pi/3}(0.14343)\right]\right\},
\]
\[
\left.I_1^k\right|_{(0,-1)}\approx -11.17447.
\]

Now we compute $I_1^k$ for $(y,x)\to(0,0)$:
\BE
\left.
I_1^k\right|_{(0,0)}=Re\left\{\eh i4k\sqrt{2}
\int^{1}_{0}\left[-e^{-i\pi/k} s^{2k}-e^{i\pi/k} s^{-2}\right]\frac{ds}{1+s^{2k}}\right\}=+\infty.
\label{genI1k}
\EE
Hence $-\left.I_1^k\right|_{(0,0)}=-\infty\,\forall\,k$. For the horizontal path depicted in Figure 13(a) we have\BE
   I_2^k=Re\int_{\alpha}\phi_2 dh=Re\int_{\alpha}{\left(\frac{i}{g}+ig\right)} dh,\label{fori2}
\EE
where $\alpha(t)=e^{it}$ with $0\leq t\leq\pi$. From (\ref{forg3}) and (\ref{fordh3}) we have
\[
   I_2^k= Re\eh i\int_{0}^{\pi}\biggl\{\biggl[-e^{\frac{-i\pi}{k}} e^{-\frac{i(k+1)t}{4k}}
   \left(\frac{y-e^{it}}{y-e^{-it}}\right)^{-\frac{k-1}{4k}}
   \left(\frac{x-e^{-it}}{x-e^{it}}\right)^{-\frac{1}{2}}\biggl.\biggl.
\]
\[
   \biggl.-e^{\frac{i\pi}{k}}e^{\frac{i(k+1)t}{4k}}\left(\frac{y-e^{it}}{y-e^{-it}}\right)^{\frac{k-1}{4k}}
   \left(\frac{x-e^{-it}}{x-e^{it}}\right)^{\frac{1}{2}}\biggl]
\] 
\BE
\biggl.\frac{i/y}{\sqrt{1+\frac{(1-y)\cos (t/2)}{\sqrt{y}(1/y+y-2\cos t)^{1/2}}}}
\frac{dt}{2\cos t-(1/y+y)}\biggl\}.
\label{i2}
\EE

We begin with $I_2^k$ and make $(y,x)\to(0,-1)$. This results in
\BE
\left.
I_2^k\right|_{(0,-1)}=-2Re\biggl\{
\int^\pi_0\cos\left(\frac{i(4\pi+t(k-1))}{4k}\right)\frac{dt}{\sqrt{1+\cos(t/2)}}\biggl\}.
\label{genI2k}
\EE

For $k=3$, (\ref{genI2k}) becomes 
\BE
\left.
I_2^3\right|_{(0,-1)}=-2Re\biggl\{
\int^\pi_0\cos\left(\frac{2\pi+t}{6}\right)\frac{dt}{\sqrt{1+\cos(t/2)}}\biggl\}.
\label{genI23}
\EE

Since $0\leq t\leq\pi$, (\ref{genI23}) rewrites as
\[
   -2\leq-2\cos\left(\frac{2\pi+t}{6}\right)\frac{1}{\sqrt{1+\cos (t/2)}}\leq2,
\]
whence 
\[
   -2\pi\leq\left.I_2^3\right|_{(0,-1)}\leq2\pi.
\]

Finally, we analyse $I_2^k$ when $(y,x)\to(0,0)$. From (\ref{i2}) it follows that
\[
 \left.
 I_2^k\right|_{(0,0)}=
-Re\eh 2\int^{\pi}_{0}\cos\left(\frac{i(4\pi-t(k+1))}{4k}\right)\frac{dt}{\sqrt{1+\cos (t/2)}}.
 \label{genIk2}
\]

For $k=3$, (\ref{genIk2}) rewrites as
\[
 \left.
 I_2^3\right|_{(0,0)}=
-Re\eh 2\int^{\pi}_{0}\cos\left(\frac{\pi-t}{3}\right)\frac{dt}{\sqrt{1+\cos (t/2)}}.
\]

Since $0\leq t\leq\pi$, then
\[
-2\leq-2\cos \left(\frac{\pi-t}{3}\right)\frac{1}{\sqrt{1+\cos (t/2)}}\leq2
\]
whence
\[
-2\pi\leq\left.I_2^3\right|_{(0,0)}\leq 2\pi.
\]
Therefore,
\[
\left.
I_1^3\right|_{(0,-1)}\approx -11.1747,
\]
\[
-2\pi\leq\left.I_2^3\right|_{(0,-1)}\leq2\pi,
\]
and then $-I_1^3>I_2^3$ at $(0,-1)$. Now
\[
-\left.I_1^3\right|_{(0,0)}=-\infty,
\]
\[
-2\pi\leq\left.I_2^3\right|_{(0,0)}\leq2\pi,
\]
whence $-I_1^3<I_2^3$ at $(0,0)$. By the {\it Intermediate Value Theorem}, there exists a point $(y^*,x^*)$ at which $-I_1^3=I_2^3$. We are ready to prove the following result, which concludes this section:
\begin{lem}
For any natural $k\ge3$ there exists a point at which $-I_1^k=I_2^k$.\label{lmm}
\end{lem}
\ \\
{\it Proof:} From (\ref{genIk1}), $\left.-I_1^k\right|_{(0,-1)}$ will be increasing with $k$ exactly when $\pi\cdot\frac{\sin(\pi/k)}{(\pi/k)}\cdot\left(-s^{-2}+\frac{1+s^{-2}}{1+s^{2k}}\right)$ is increasing with $k$. But this is obvious to the function $\frac{\sin(\pi/k)}{(\pi/k)}$, and also to the function $-s^{-2}+\frac{1+s^{-2}}{1+s^{2k}}$. Since both are positive, then $\left.-I_1^k\right|_{(0,-1)}$ is increasing with $k$. Moreover, for all $k\ge3$ we have
\[
   -2\pi\le I_2^k|_{(0,-1)}\le 2\pi,
\]
\[
   -2\pi\le I_2^k|_{(0,0)}\le 2\pi,
\]
because these inequalities hold for $k=3$, and the same computations lead to the general case. Hence $-I_1^k|_{(0,-1)}\ge-I_1^3|_{(0,-1)}>I_2^k|_{(0,-1)}$, and also $-I_1^k|_{(0,0)}=-\infty<I_2^k|_{(0,0)}$.\hfill q.e.d.  

\section{Embeddedness}

Now we have a compact Riemann surface $\ovl{S}$ given by (\ref{fg}), where $(x,y)$ is the point that Lemma \ref{lmm} refers to. If $B\subset\ovl{S}$ denotes the branch points of $g$, then the ends of $\ovl{S}$ are $\{p_1,\dots,p_{2k}\}=z^{-1}(1)\setminus B$, as shown in the next subsection. We have $dh$ defined on $\ovl{S}$, with a {\it local} expression given by (\ref{eqdh}). From Theorem \ref{WR}, this defines a minimal immersion $X:R\to{\mathbb{E}}$, where $R=\ovl{S}\setminus\{p_1,\dots,p_{2k}\}$. The purpose of this section is to prove that $X$ is an embedding.

\subsection{Poles and zeros of $dg$}
\label{secpolos}

We know that $deg(dg)=-\chi_{\bar{S}}=4k-2$. According to Figure 4, $dg$ should have exactly $k+(2+2)k=5k$ poles. Notice that the multiplicity is always included in our analysis. 

Since $deg(dg)=$ Nr zeros$(g)-$ Nr poles$(g)$, then $dg$ has exactly $9k-2$ zeros. Now we are going to locate these zeros by geometric arguments. They can be checked analytically from (\ref{fg}).

According to Figure 4, $k-2$ of the zeros are at the saddle point marked with $\times$ in Figure 2 where $z=y$. We have $k$ 4-fold saddle points along each horizontal closed curve of symmetry, which gives a subtotal of $2k$ points. Each of them adds $2$ zeros for $dg$. Hitherto we have $5k-2$ zeros. 

Somewhere along each ray that departs from the $\times$-point, the unitary normal has an inflexion, hence $dg=0$ there. By counting all such rays, at the top and at the bottom of the fundamental piece, we arrive at the $4k$ remaining zeros that finally totalize $9k-2$.

\subsection{The U-curve}
\label{ucurve}
Soon we shall see that the $U$-curve is not of symmetry. However, in Section 4 we explained that our $dh$ is a {\it particular} case of a more general formula that should involve $z$, $dz$, $g$ and $dg$. Therefore, it is highly possible that, for each $k$, the surfaces $ST_{2k}$ are representatives of a continuous {\it two}-parameter family of surfaces. This family is very likely to have more symmetric members for which the $U$-curve {\it is} of symmetry.

Therefore, herewith we present arguments under this assumption. It will be slacked later, but the reader will notice how simple the arguments are in this case. Most of these arguments remain valid in the slacked case, and so we ease the understanding of our proof. 

We named $P$ the fundamental piece of $ST_{2k}$. In Figure 13, the shaded region represents a fundamental domain of $P$.

\begin{figure}[!htb]
\begin{minipage}[b]{0.4\linewidth}
\includegraphics[width=\linewidth]{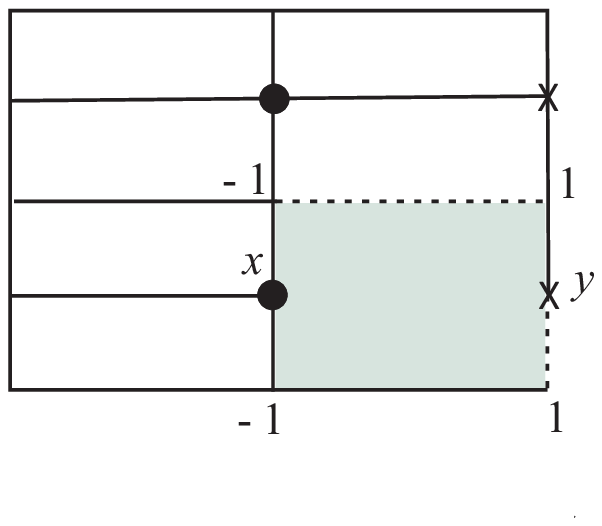}
\caption{A fundamental domain in the torus $T$.}
\end{minipage}\hfill
\begin{minipage}[b]{0.40\linewidth}
\includegraphics[width=\linewidth]{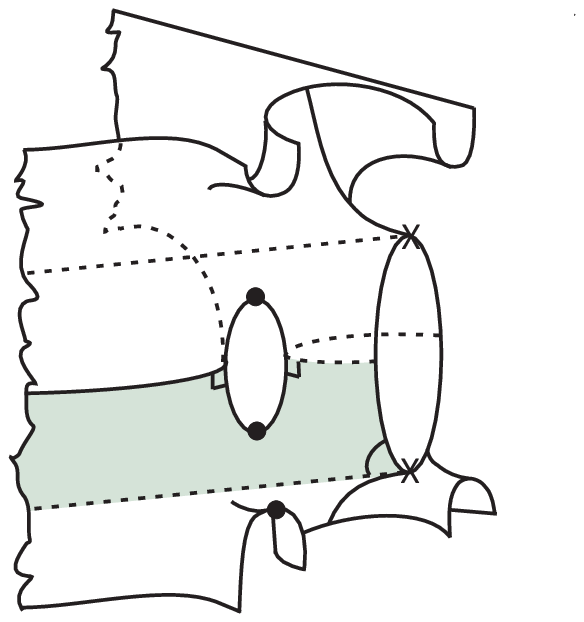}
\caption{A fundamental domain of $P$.}
\end{minipage}
\end{figure}

Let ${\cal{R}}$ be the circle in Figure 14 left. The image ${\cal{D}}=g({\cal{R}})$ is depicted in Figure 14 right. It is contained in a hemisphere of $S^2=\hC$. Hence, there is a direction in which the orthogonal projection of $X({\cal{R}})$ is an immersion. For instance, direction $Ox_2$. This way $(x_1,x_3):{\cal{R}}\to\R^2$ is an immersion when restricted to the interior of ${\cal{R}}$. The image $(x_1,x_3)({\cal{R}})\subset\R^2$ has one out of four basic features described in Figure 15.

\begin{figure}[h!]
\centering
\includegraphics[scale=0.80]{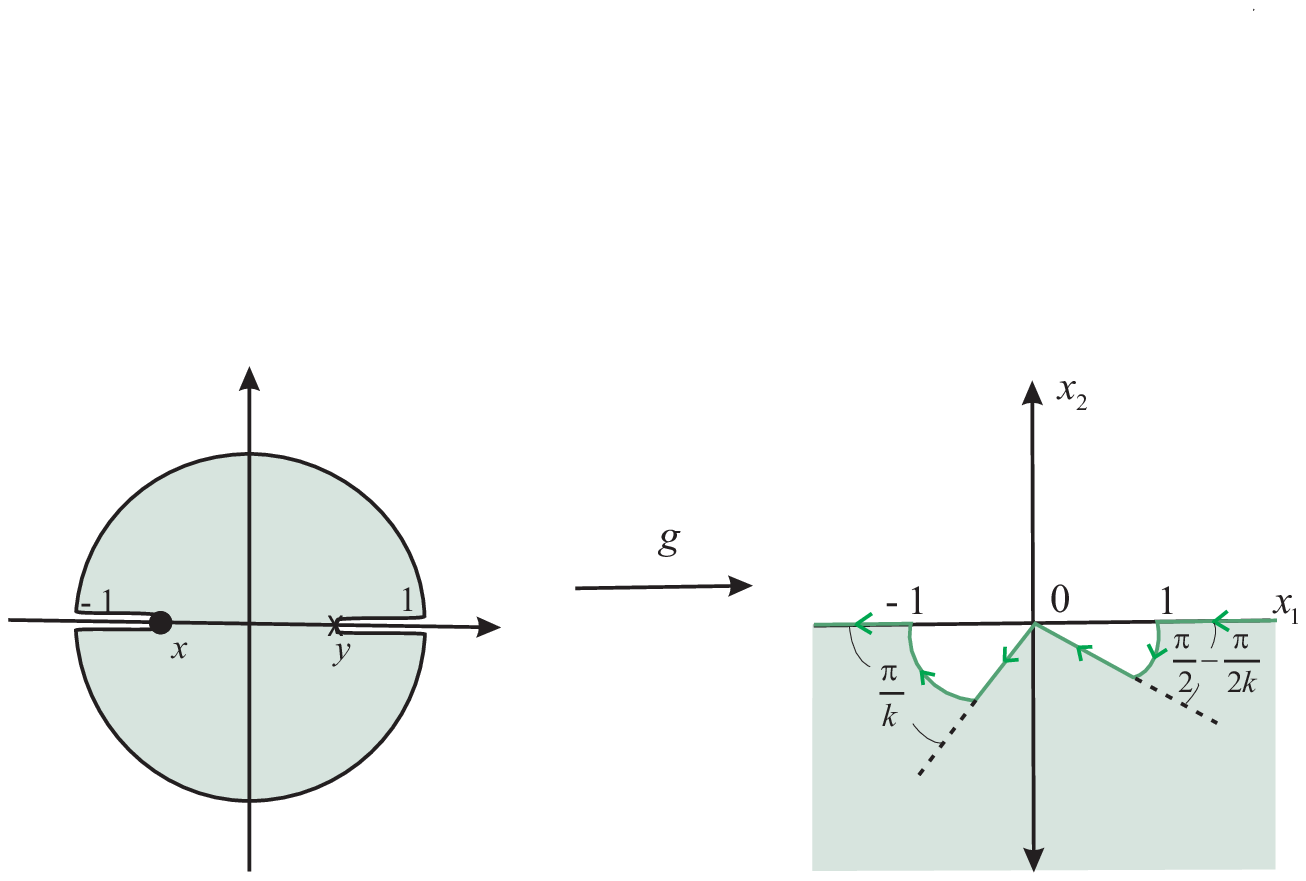}
\caption{Image ${{\cal{D}}=g({\cal{R}})}$.}
\end{figure} 

\begin{figure}[h!]
\centering
\includegraphics[scale=0.80]{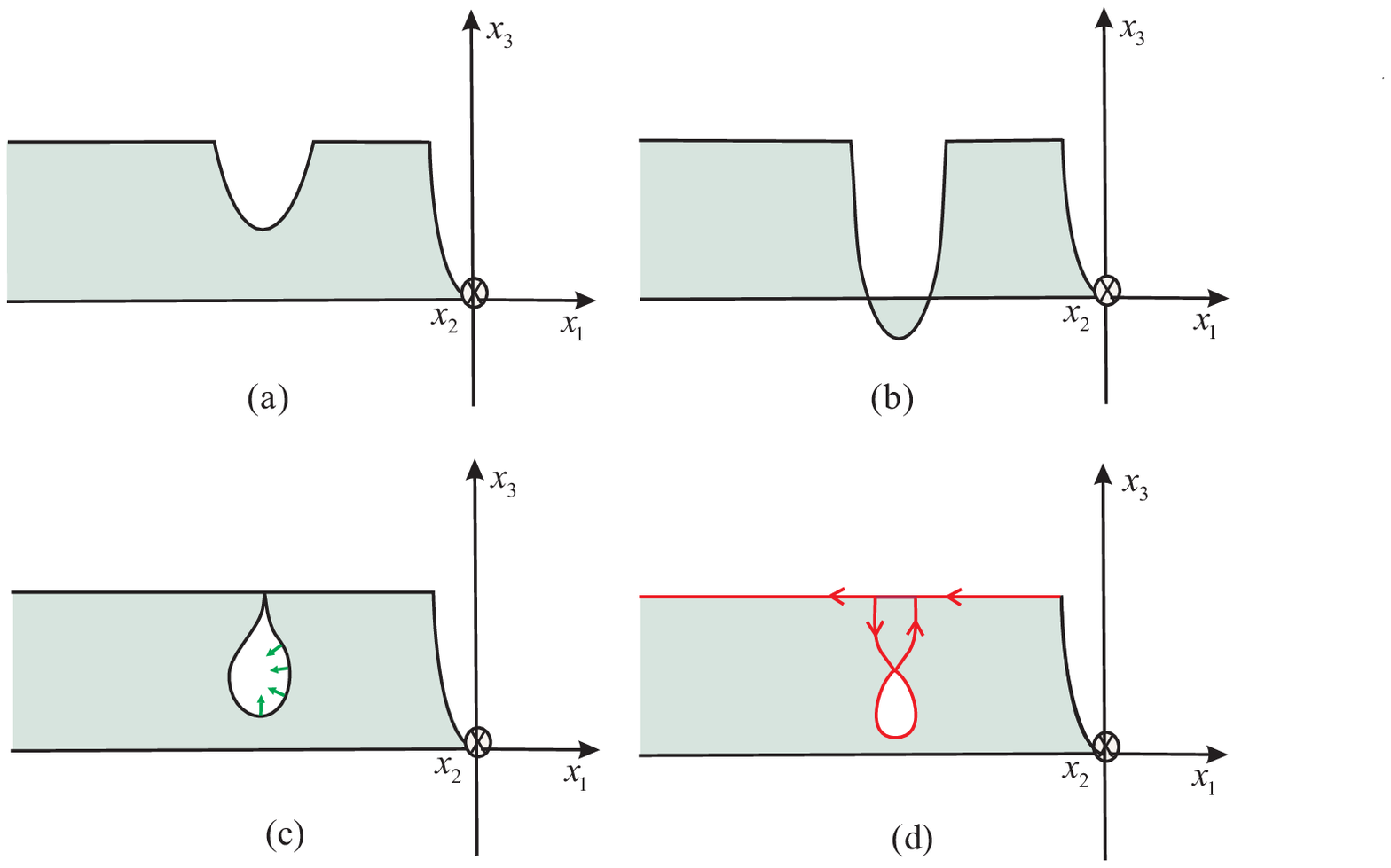}
\caption{Four basic features of $(x_1,x_3)({\cal{R}})$.}
\end{figure} 

If the $U$-curve is planar, then it is convex. Indeed, the Gaussian curvature 
\[
   K=-\biggl(\frac{2}{|g|+1/|g|}\biggl)^4\biggl|\frac{dg/g}{dh}\biggl|^2
\]
vanishes on $U$ exactly where $dg=0$ there. But in Subsection~\ref{secpolos} we saw that this never happens. 

Since $(x_1,x_3)|_{\cal{R}}$ is an immersion, it is open and continuous. Hence, $(x_1,x_3)({\cal{R}})$ is an open connected subset of $\R^2$, which discards Figure 15(b). Figures 15(c) and 15(d) are also discarded because $g$ is injective along $U$. We finally remain with Figure 15(a).

Let $G$ be the interior of $(x_1,x_3)({\cal{R}})$. From Theorem \ref{simcon}, $G$ is simply connected. Of course, ${\cal{R}}\subset\hC$ and $(x_1,x_3)$ extends continuously to $(x_1,x_3):{\bar{\cal{R}}}\to\hC$ by taking $(x_1,x_3)(\infty)=\infty$. The pre-image of any point in $(x_1,x_3)({\cal{R}})$ is a finite subset of ${\cal{R}}$, otherwise it would have an accumulation point at the boundary of ${\cal{R}}$, but $\partial G$ consists of monotone curves. Hence, $(x_1,x_3)|_{\cal{R}}$ is a covering map of the simply connected $G$. 

Namely, $(x_1,x_3)|_{\cal{R}}$ is injective. Therefore, $(x_1,x_2,x_3):{\cal{R}}\to\R^3$ is a graph. 

\begin{rmk}
\label{obsUns}
The $U$-curve is not a symmetry curve. Otherwise, it would be in $Ox_1x_3$ and the conjugate minimal surface would have a straight segment perpendicular to this plane. But $idh$ is the 3rd coordinate of the conjugate, and so its real part should be zero along the segment. However, in Section~\ref{secdh} we saw that $dh\not\in\R\cup i\R$ along this curve.
\end{rmk}

\subsection{Embeddedness proof}

We shall some ideas from \cite{Valerio2}. Two copies of the fundamental domain are represented in Figure 16, together with corresponding $z$-image, namely two superimposed unitary disks.

\begin{figure}[h!]
\centering
\includegraphics[scale=0.80]{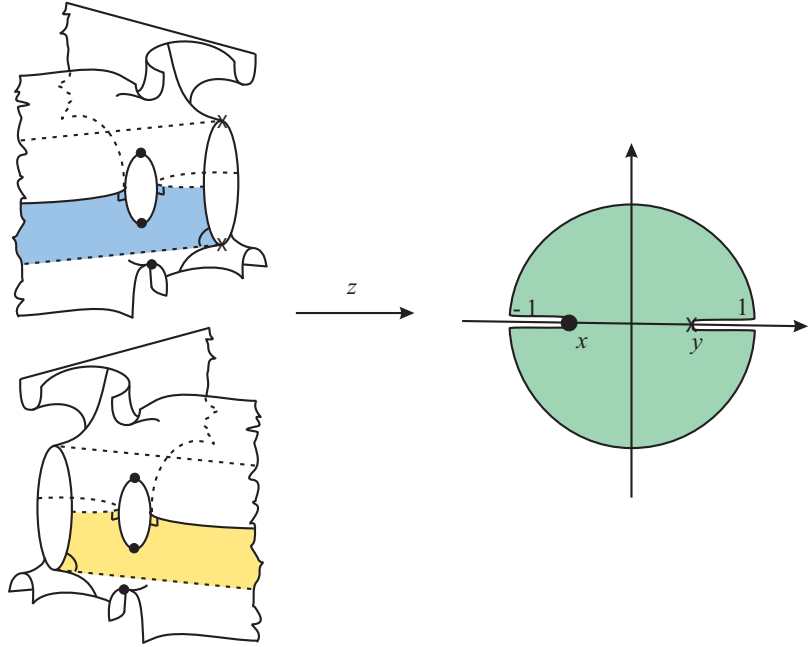}
\caption{Two copies of the fundamental domain and their $z$-image.}
\end{figure} 
 
The corresponding $g$-image is depicted in Figure 17. Let $\Gamma$ be the $z$-image of the $U$-curve. Then $\Gamma$ separates the superimposed disks into two disjoint components $\A$ and $\B$.

\begin{figure}[ht!]
\centering  
\includegraphics[scale=0.80]{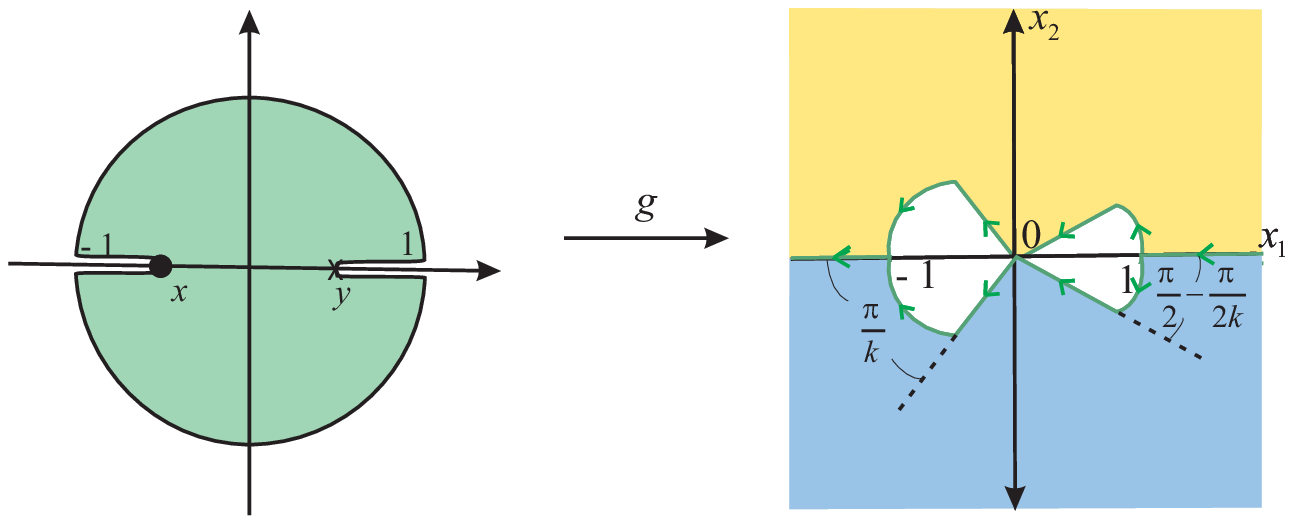}
\caption{The $g$-image of $\A\cup\B$.}
\end{figure} 

Figure 17 shows that $g(\A)$ is the conjugate of $g(\B)$, whereas $g(\Gamma)=\hat{\R}\setminus(-1,1)$. 

In Subsection~\ref{ucurve} we saw that $(x_1,x_3):\A\to\R^2$ and $(x_1,x_3):\B\to\R^2$ are both immersions. Again, Figure 15 shows the four possible features of $(x_1,x_3)(\A)$. Since $dg$ never vanishes along the $U$-curve, then Figures 15(c) and 15(d) cannot occur. We also discard Figure 15(b) by the same arguments presented in Subsection~\ref{ucurve}. So, there remains Figure 15(a) and we conclude that $(x_1,x_2,x_3):\A\to\R^3$ and $(x_1,x_2,x_3):\B\to\R^3$ are both graphs.

Except for their boundaries, none of the graphs can intercept the other. Otherwise, they would either be tangent, or we could make them tangent by displacing both graphs in opposite directions along $Ox_2$. This would contradict the {\it Maximum Principle for Minimal Surfaces}. 

Therefore, both graphs intersect only at their coinciding $U$-curves, and what we have is an embedded double-piece of minimal surface contained in a wedge of $\R^3$. Its angle is $2\pi/k$, and except for the $U$-curve, the double-piece has its rays and reflectional symmetry curves on the faces of the wedge. The whole $ST_{2k}$ is then generated by successive rotations around the rays and reflections on the symmetry curves. Therefore, it has no self-intersections. Since $X$ is proper, then it is an embedding.

{\rm A.J. Yucra Hancco\\ DM, UFSCar, \\ Rua Washington Lu{\'{\i}}s km 235,\\13565-905 S\~ao Carlos, SP, Brazil\\ E-mail: {\it alvaro@dm.ufscar.br}}\\
\\
{\rm G. A. Lobos\\ DM, UFSCar, \\ Rua Washington Lu{\'{\i}}s km 235,\\13565-905 S\~ao Carlos, SP, Brazil\\ E-mail: {\it lobos@dm.ufscar.br}}\\
\\
{\rm V. Ramos Batista\\ CMCC, UFABC, \\Rua Santa Ad\'elia 166, Bl.A-2,\\09210-170 Santo Andr\'e, SP, Brazil\\ E-mail: {\it valerio.batista@ufabc.edu.br}}


\begin{thebibliography}{15}
\itemsep = 0.0 pc
\parsep  = 0.0 pc
\parskip = 0.0 pc    
\bibitem{Conway} J.B. Conway, {\it Function of one complex variable}. New York. Springer-Verlag, (1979).
\bibitem{Costa} C.J. Costa, {\it Example of a complete minimal immersion  in $\R^3$ of genus one and three embedded ends}. Bol. Soc. Brasil. Mat. {\bf 15}, 47-54 (1984).
\bibitem{Foster} O. Foster, {\it Lectures on Riemann Surfaces}.  New York. Springer-Verlag, (1981).
\bibitem{HMM} L. Hauswirth, F. Morabito, M. Rodriguez. {\it An end-to-end-construction for singly periodic minimal surfaces}. Pacific J. Math. 241, 1--61 (2009).
\bibitem{Hoffman1} D. Hoffman, W.H. Meeks, {\it Embedded minimal surfaces of finite topology}. Ann. Math. {\bf 131}, 1-34 (1990).
\bibitem{Hoffman} D. Hoffman, H. Karcher, {\it Complete embedded minimal surfaces of finite total curvature}. Encyclopaedia Math. Sci. {\bf 90} (1997).
\bibitem{JorgeM} L.P.M. Jorge, W.H. Meeks, {\it The topology of complete minimal surfaces of finite total Gaussian curvature}. Topology {\bf 2}, 203-221 (1983).
\bibitem{Karcher1} H. Karcher, {\it Embedded minimal surfaces derived from Scherk's examples}. Manuscripta Math. {\bf 62}, 83-114, Bonn, (1988).
\bibitem{Karcher} H. Karcher, {\it Construction of minimal surfaces, Surveys in Geometry}. University of Tokyo. Lecture Notes, vol. 12, pp 1-96, SFB256, Bonn, (1989).
\bibitem{Lopez} F.J. L\'opez, F. Mart\'\i n, {\it Complete minimal surfaces in $\R^3$}. Publ. Mat. {\bf 43}, 341-449 (1999).
\bibitem{Lopez2} F.J. L\'opez, A. Ros, {\it On embedded complete minimal surfaces of genus zero}. J. Differ. Geom. {\bf 33}, 293-300 (1991).
\bibitem{Martin} F. Martin, V. Ramos Batista, {\it The embedded singly periodic Scherk-Costa surfaces}. Math. Ann. {\bf 336}, 155-189 (2006).
\bibitem{Meeks} W.H. Meeks, M. Wolf, {\it Minimal surfaces with the area growth of two planes; the case of infinite symmetry}. J. Amer. Math. Soc. {\bf 20}(2), 441-465 (2007).
\bibitem{Nitsche} J.C.C. Nitsche, {\it Lectures on minimal surfaces}. Cambridge University Press, Cambridge, (1989).
\bibitem{Osserman} R. Osserman, {\it A Survey of Minimal Sufaces}. 2nd edn. Dover, New York (1969).
\bibitem{Valerio1} V. Ramos Batista, {\it A family of triply periodic Costa surfaces}. Pacific J. Math. {\bf 212}, 347-379 (2003).
\bibitem{Schoen1} A.H. Schoen, {\it Infinite periodic minimal surfaces without selfintersections}. NASA technical note, No. D-5541 (1970).
\bibitem{Schoen} R. Schoen, {\it Uniqueness, symmetry, embeddedness of minimal surfaces}.  J. Differ. Geom. {\bf 18}, 791-809 (1983).
\bibitem{Valerio2} M. F. da Silva, V. Ramos Batista, {\it Scherk saddle towers of genus two in $\R^3$}. Geom. Dedicata {\bf 149}, 59-71 (2010).
\bibitem{Traizet} M. Traizet, {\it Construction de surfaces minimales en recollant des surfaces de Scherk}, Ann. Inst. Fourier {\bf 46}, 1385-1442 (1996).
\end{thebibliography}
\end{document}